\documentclass[twoside,11pt]{article}
\usepackage{graphicx}
\usepackage{amssymb}
\usepackage{amsthm}
\newfont{\bbb}{msbm10 scaled\magstep 1}

\newcommand{\conv}{{\rm conv}}

\newcommand{\bd}{{\rm bd}}
\newcommand{\width}{{\rm width}}

\newtheorem{lem}{Lemma}

\newtheorem{cla}{Claim}

\title{\bf  Approximation of Spherical Bodies of Constant Width and Reduced Bodies}
\date{}
\begin{document}\maketitle

\vskip -0.5cm 
\centerline
{\bf Marek Lassak}

\pagestyle{myheadings} \markboth{ M. Lassak / Approximation of reduced bodies on sphere}{ M. Lassak / Approximation of reduced bodies on sphere}

\vskip 0.9cm 
 
\baselineskip = 16pt

\noindent 
{\bf Abstract.} We present a spherical version of the theorem of Blaschke that every body of constant width $w < \frac{\pi}{2}$ can be approximated by a body of constant width $w$ whose boundary consists only of pieces of circles of radius $w$ as well as we wish in the sense of the Hausdorff distance.
This is a special case of our theorem about approximation of spherical reduced bodies.

\baselineskip = 16pt



\section{Introduction}

\noindent
A theorem of Blaschke says that for every convex body of constant width $w$ in the Euclidean plane $E^2$ and every $\varepsilon >0$ there exists a convex body of constant width $w$ whose boundary consists only of pieces of circles of radius $w$ such that the Hausdorff distance between the two bodies is at most $\varepsilon$ (see \cite{Bl} and also \S 65 of \cite {BF}). 
A generalization of this fact for normed planes is given in \cite{L1}, where also approximation of reduced bodies is considered.
Corollary at the end of this note presents an analog of this theorem for bodies of constant width on the sphere $S^2$, while Theorem gives a general version for reduced convex bodies on $S^2$.

The method of the proof of Theorem is similar to the proofs from \cite{Bl} (by Blaschke) and \cite{L1}.
In order to facilitate the reader a comparison of the present proof with this from  \cite{L1} for normed planes, in Section 4 we use the notation from \cite{L1}.
In the present paper we need a number of lemmas and claims for the spherical situation.
They are given in Sections 2 and 3.

\section{Auxiliary facts from spherical geometry}

Let $S^2$ be the unit sphere of the three-dimensional Euclidean space $E^3$. 
By a {\it great circle} we mean the intersection of $S^2$ with any 
two-dimensional subspace of $E^{3}$. 
By a pair of {\it antipodes} we mean any pair of points obtained as the intersection of $S^2$ with a one-dimensional subspace of $E^3$.
Observe that if two different points $a, b$ are not antipodes, then there is exactly one great circle containing them.
By the {\it arc of great circle}, or shortly {\it arc}, $ab$ connecting them we understand the shorter part of the great circle between  $a$ and $b$. 
By the {\it distance} $|ab|$, of these points we mean the length of $ab$. 

The set of points of $S^2$ whose distance from a point $c \in S^2$ equals (is at most) $r \leq \frac{\pi}{2}$ is called the {\it circle} (respectively, {\it disk}) {\it of radius $r$ and center $c$}.
Disks of radius $\frac{\pi}{2}$ are called {\it hemispheres}.
If $p$ belongs to a circle, then the set of points which are at distances  at most $\frac{\pi}{2}$ from $p$ is called a {\it semicircle}.
We say that $p$ is the {\it center} of this semicircle.

A subset of $S^2$ is called {\it convex} if it does not contain any pair of antipodes and if together with every two points it contains the arc connecting them.
For a set $A$ contained in the interior of a hemisphere we denote by $\conv (A)$ the smallest convex set containing $A$.
By {\it a convex body} we mean a closed convex set with non-empty interior. 
Let $p$ be a boundary point of a set $A$ contained in the interior of a hemisphere. 
We say that a hemisphere $K$ containing $A$ {\it supports} it if $A \cap \bd (K)$ in non-empty.  
If $p \in  A \cap K$, we say that $K$ supports $A$ at $p$.

\begin{lem} 
Let $A \subset S^2$ be a closed set with non-empty interior.
If at every boundary point the set $A$ is supported by a hemisphere, then $A$ is convex. 
\end{lem} 
 
\begin{proof}
The proof is analogous to the proof of Theorem 9 of \cite{E}.
Namely, suppose that $A$ is not convex.
Then there are points $x_1, x_2 \in A$ and a point $y \in x_1x_2$ such that $y \not \in A$. 
There exists an interior point $z$ of $A$ such that $z \not \in x_1x_2$.
%There exists a point $z \in inter (A)$ such that $z \not \in x_1x_2$.
%There exists a point $z \in \inter (A)$ such that $z \not \in x_1x_2$.
Hence there is $x_f \in \bd (A) \cap yz$.
Every great circle through $x_f$ separates one of the points $x_1, x_2, z$ from at least one other of these points.
So there is no supporting hemisphere of $A$ through $x_f$.
A contradiction.
\end{proof}

Recall a few notions and facts from \cite{L2}. 
We say that $e$ is an \emph{extreme point} of $C$ provided the set $C \setminus \{e\}$ is convex.
If hemispheres $G$ and $H$ are different and their centers are not antipodes, then $L = G \cap H$ is called a {\it lune}. 
The semicircles bounding $L$ and contained respectively in $\bd (G)$ and $\bd (H)$ are denoted by $G/H$ and $H/G$. 
The {\it thickness} $\Delta (L)$ of $L$ is defined as the distance of the centers of $G/H$ and $H/G$.  
For every hemisphere $K$ supporting a convex body $C$ we find hemispheres $K^*$ supporting $C$ such that the lunes $K\cap K^*$ are of the minimum thickness.
At least one such a $K^*$ exists.
By {\it the width ${\rm width}_K (C)$ of $C$ determined by} $K$ we mean the thickness of the lune $K \cap K^*$. 
It changes continuously, as the position of $K$ changes. 
From Part I of Theorem 1 of \cite{L2} we know that such $K^*$ is unique if ${\rm width}_K (C) < \frac{\pi}{2}$.

When we go on $\bd (C)$, then always counterclockwise.
Let us recall, but in a slightly different form than in \cite{LaMu}, that if hemispheres $X, Y, Z$ support a convex body $C$ in this order, then we write $\prec \!XYZ$. 
Let $X$ and $Z$ be hemispheres supporting $C$ at $p$ and let $\preceq \! XYZ$ for every hemisphere $Y$ supporting $C$ at $p$.
Then $X$ is said to be the {\it right supporting} hemisphere at $p$
and $Z$ is said to be the {\it left supporting} hemisphere at $p$. 
 
By the {\it thickness} $\Delta (C)$ of $C$ we mean the minimum of $\width_K (C)$ over all supporting hemispheres $K$ of $C$. 

Assume that a lune $L$ of thickness $\Delta C)$ contains a convex body $C$.
Then by Claim 2 of \cite{L2} both the centers of the semicircles bounding $L$ belong to $\bd (C)$. 
We call such a lune $L$ a {\it supporting lune of} $C$. 
We say that the lune {\it supports $C$ at} the mentioned centers.
The arc connecting the centers of the semicircles bounding $L$ is called a {\it thickness chord of} $C$. 
This notion is an analog of the notion of a thickness chord of a convex body in $E^d$ considered by many authors under this name or without any name, as for instance in \cite{BF}. 

If for all hemispheres $K$ supporting $C$ the numbers ${\rm width}_K (C)$ are equal, we say that $C$ is {\it of constant width}. 
Of course, a convex body $C \subset S^2$ with $\Delta (C) < \frac{\pi}{2}$ is of constant width if and only if for every hemisphere $K$ supporting $L$ the lune $K \cap K^*$ has thickness $\Delta (C)$.

\begin{lem} 
Consider a non-degenerate spherical triangle $abd \subset S^2$ and a point $e \in ab$ such that $de$ is orthogonal to $ab$.
The distance between $d$ and $ab$ is at most 

$${\rm arc sin}\frac{\tan \frac{1}{2} |ab|}{{\tan \frac{1}{2} \angle adb}}\; .$$
\end{lem}

\begin{proof} 
Let $e'$ be the center of $ab$ and $d'$ be such that $|d'e'| = |de|$ with $d'e'$ orthogonal to $ab$.
Clearly, $\angle e'd'a = \frac{1}{2}\angle ad'b$ and $|ae'| = \frac{1}{2}|ab|$.
Since $ae'd'$ is a right triangle, $\sin |d'e'| \cdot \tan \frac{1}{2} \angle adb = \tan \frac{1}{2} |ab|$ and our thesis for $d'$ instead of $d$ is true.
By $|d'e'| = |de|$ we get the thesis for our $d$.
\end{proof}

Here is a spherical version of the classic Blaschke selection theorem.

\begin{cla} From every sequence of convex bodies on $S^d$ of thickness at most a fixed constant smaller than $\pi$ we may select a subsequence which tends to a spherical convex body.
\end{cla}

\begin{proof}
First we select a subsequence which is contained in a certain spherical disk of a radius below $\pi$. 
Then from this subsequence we select the final subsequence (see \cite{L2} p. 558 and \cite{HA}).
\end{proof}

\section{Reduced bodies on sphere}

After \cite{L2} we say that a spherical convex body $R$ is {\it reduced} if $\Delta (Z) < \Delta (R)$ for every convex body $Z \subset R$ different from $R$. 
For basic properties of reduced bodies on $S^2$ see \cite{LaMu}.
The class of spherical reduced bodies is larger than the class of bodies of constant width.
Earlier many papers considered reduced bodies in the Euclidean and normed $d$-dimensional spaces; see the survey articles \cite{LasMart1} and \cite{LasMart2}.

Later tacitly assume that all considered reduced bodies $R \subset S^2$ are of thickness below $\frac{\pi}{2}$.

Theorem 3.1 of \cite{LaMu} presents the boundary structure of a reduced body $R \subset S^2$. 
Namely, {\it assume that $M_1, M_2$ are supporting hemispheres of $R$ with $\width_{M_1}(R) = \Delta(R) = \width_{M_2}(R)$ such that $\width_M(R) > \Delta(R)$ for every $M$ fulfilling $\prec M_1MM_2$.
Consider the lunes $L_1 = M_1 \cap M_1^*$ and $L_2 = M_2 \cap M_2^*$.
Then the arcs $a_1a_2$ and $b_1b_2$ are in $\bd(R)$, where $a_i$ is the center of $M_i/M_i^*$ and $b_i$ is the center of $M_i^*/M_i$ for $i=1,2$.}

Denote by $c$ the intersection point of the arcs $a_1a_2$ and $b_1b_2$. 
The union of the triangles $a_1a_2c$ and $b_1b_2c$ is called a {\it butterfly}, while $a_1a_2$ and $b_1b_2$ its {\it arms}.

Consider any maximum piece $\buildrel \frown \over {gh}$ of the boundary of $R$ which does not contain any arc 
Let $K$ be any hemisphere supporting $R$ at a point of $\buildrel \frown \over {gh}$.
But at $g$ we agree only for the right, and at $h$ we agree only for the left.
By Part I of Theorem 1 of \cite{L2} for every $K$ there exists exactly one hemisphere $K^*$ supporting $R$ such that the lune $K \cap K^*$ has thickness $width_K(R)$.
From this theorem we also conclude that $K^*$ touches $R$ at a unique point. 
In particular, for the right hemisphere supporting $R$ at $g$ denote this unique point by $g'$, and for the left hemisphere supporting $R$ at $h$ by $h'$.

\begin{cla}
For any hemisphere $K$ supporting $R$ at a point of $\buildrel \frown \over {gh}$ we have 
$\width_K (R) = \Delta(R)$.
\end{cla}

\begin{proof}
Clearly, $\width_K (R)$ cannot be smaller than $\Delta (R)$, because then we get a contradiction with the definition of the thickness of $R$.

Suppose that $\width_K (R) > \Delta(R)$. 
Then by continuity arguments (see Theorem 2 of \cite{L2}) we obtain that if $K$ supports $R$ at a point different from $g$ and $h$, then there are supporting hemispheres $K_1$ and $K_2$ with $\prec K_1KK_2$ such that for every supporting hemisphere $H$ fulfilling $\prec K_1HK_2$ we have $\Delta(H) > \Delta(R)$.
By the just recalled Theorem 3.1 of \cite{LaMu}, the piece $\buildrel \frown \over {gh}$ of the boundary of $R$ contains an arc.
A contradiction with the assumption on $\buildrel \frown \over {gh}$ from the paragraph preceding our claim.
If $K$ supports $R$ at $g$ or $h$, then similarly we show that $\buildrel \frown \over {gh}$ contains an arc whose one end-point is $g$ or $h$, which again gives a contradiction with the choice of $\buildrel \frown \over {gh}$.

By the two preceding paragraphs we get $\width_K (R) = \Delta(R)$.
\end{proof}

From this claim and Lemma 2.2 of \cite{LaMu} we see that all the points at which all our hemispheres $K^*$ touch $R$ form a  curve $\buildrel \frown \over {g'h'}$ being a piece of $\bd(R)$.
We call it the curve {\it opposite to the curve} $\buildrel \frown \over {gh}$.
Vice-versa, from this lemma we obtain that $\buildrel \frown \over {g'h'}$ determines $\buildrel \frown \over {gh}$.
So we say that $\buildrel \frown \over {gh}$ 
and $\buildrel \frown \over {g'h'}$ is a {\it pair of opposite curves of constant width $\Delta(R)$}. 
Let us summarize the above consideration as the following claim analogous to Corollary 11 of \cite{FaLa} on the situation in a normed plane.

\begin{cla}
The boundary of a reduced body $R \subset S^2$ consists of countably many pairs of arms of  butterflies and of countably many pairs of opposite pairs of curves of constant width $\Delta(R)$.
\end{cla}

Clearly any particular pair of opposite pieces of curves of constant width $\Delta(R)$ consists of end-points of thickness chords; they are the centers of the semicircles bounding the lunes of thickness $\Delta(R)$ supporting $R$ at the points of these pieces of curves of constant width.

In particular, when $R$ from Claim 3 is a body of constant width, then its boundary is the union of one pair of curves of constant width $\Delta(R)$. 
We may present $\bd (R)$ on infinitely many ways as such an union.
Each time the curves end at the end-points of a thickness chord of $R$.

\section{Approximation of reduced bodies}

\vskip0.1cm
\noindent
{\bf Theorem.}
{\it Let $R \subset S^2$ be a reduced body 
of thickness below 
$\frac{\pi}{2}$.
For any $\varepsilon >0$ there exists a re- duced body $R_\varepsilon \subset S^2$ with $\Delta (R_\varepsilon) = \Delta(R)$  whose boundary consists only of arms of butterflies and arcs of circles of radius $\Delta(R)$ such that the Hausdorff distance between $R_\varepsilon$ and $R$ is at most~$\varepsilon$.}

\begin{proof}
We omit the trivial case when the boundary of $R$ consists only of arms of butterflies. 

Take any positive $\varepsilon < \frac{1}{2}\pi$. 
Put $\rho_\varepsilon = 2 \cdot {\rm arc tan} (\sin \varepsilon)$.

Consider any pair $F, G$ of opposite curves of constant width $\Delta (R)$ in the boundary of $R$.
Exceptionally, when $R$ is a body of constant width, we divide ${\rm bd}(R)$ into a pair of curves $F$ and $G$ by an arbitrary thickness chord. 
Denote the endpoints of $F$ by $f'$ and $f''$, and the endpoints of $G$ by $g'$ and $g''$, in both cases according to the positive orientation.

\vskip0.2cm
Part A. The aim of this part is to construct the set $R_\varepsilon$.

\vskip0.1cm
From Claim 3 we know that $\bd(R)$ contains countably many pairs of opposite pieces of curves 
of constant width $\Delta(R)$.
For each such a pair $F, G$ we provide a number of different thickness 
chords $f_1g_1, \dots , f_ng_n$ of $R$ such that $f_1, \dots, f_n \in F$ (with $f_1 = f'$ and $f_n = f''$), and $g_1, \dots, g_n \in G$ (with $g_1 = g'$ and $g_n = g''$), taking care that $|f_if_{i+1}|$ and  $|g_ig_{i+1}|$ be below $  \varepsilon$
for $i=1, \dots , n-1$, and that the positively oriented angle between every two successive chords is below $\frac{\pi}{2}$. 
Clearly, some of points $f_1, \dots, f_n$ (some of $g_1, \dots, g_n$) may coincide.
Let $c_i \in \Gamma_i$ be a

Denote by $o_i$ the intersection point of $f_ig_i$ and $f_{i+1}g_{i+1}$ for $i=1, \dots, n-1$  (see Figure).
Moreover, denote by $\Phi_i$ (respectively, $\Gamma_i$) the angle between the rays from $o_i$ through $f_i$ and $f_{i+1}$  (respectively, through $g_i$ and $g_{i+1}$) for $i=1, \dots, n-1$.
point of intersection of circles of

\begin{figure}[htbp]
\hskip2.6cm \includegraphics[width=11.05cm,height=7.8cm]{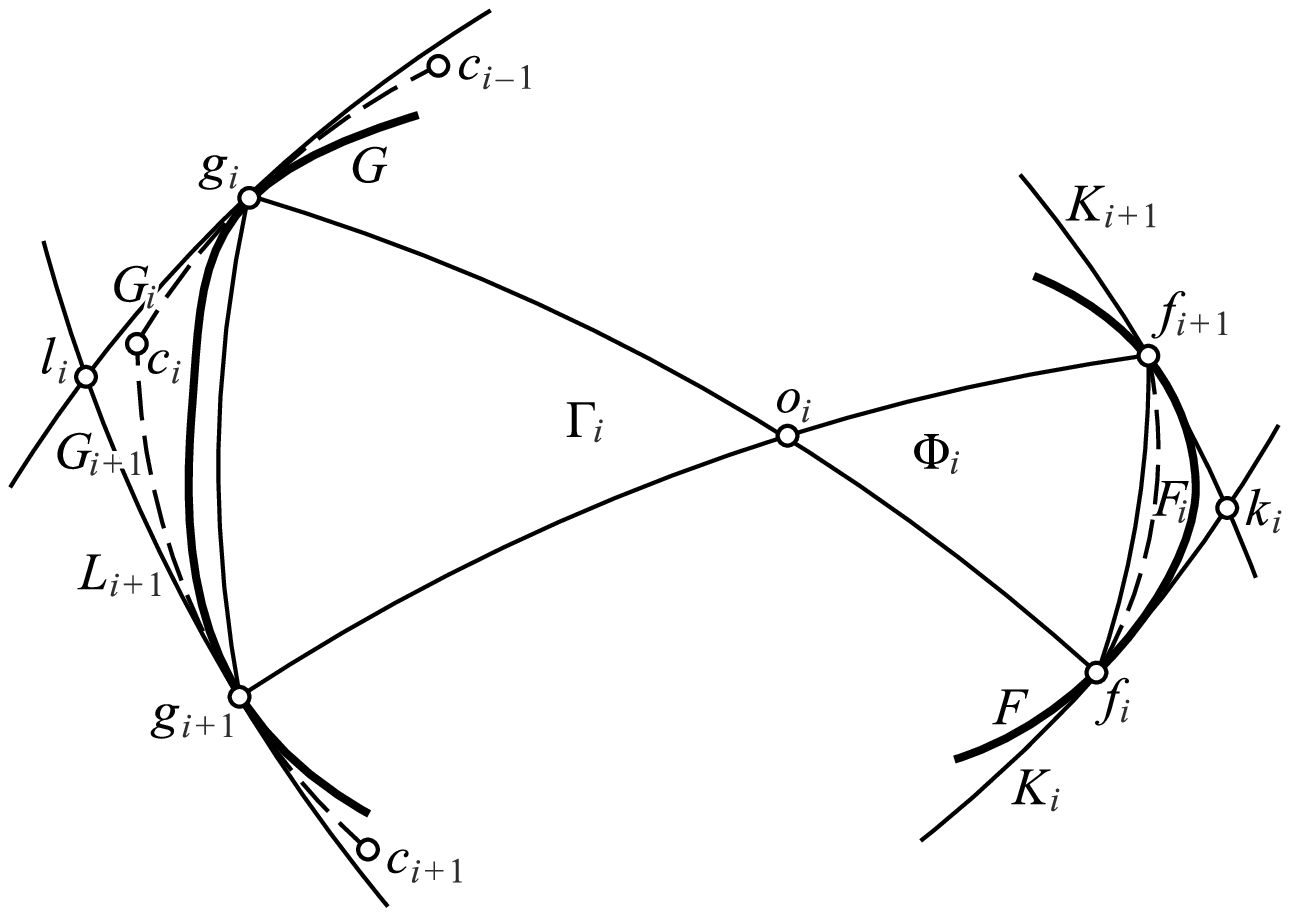} \\  
\vskip-0.1cm
\centerline{Figure: Illustration to the proof of Theorem}
\end{figure}

{\ }

\noindent
radius $\Delta(R)$ with centers $f_i$ and $f_{i+1}$ (so $c_i$ is in equal distances from $f_i$ and $f_{i+1}$).
Such $c_i$ exists since, thanks to Claim 2, we have $|f_ig_{i+1}| \leq \Delta(R)$ and $|f_{i+1}g_i| \leq \Delta(R)$.
Moreover, by $c_0$ we mean $g'$, and by $c_n$ we mean $g''$.

For every $i\in \{ 1, \dots, n-1 \}$ take the arc $F_i$ of the circle $F_i^\circ$ of radius $\Delta(R)$ with center $c_i$ and endpoints $f_i$ and $f_{i+1}$ which is in $\Phi_i$.
Moreover, for very $i\in \{ 1, \dots, n \}$ take the arc $G_i$ of the circle $G_i^\circ$ of radius $\Delta(R)$ with center $f_i$ which begins at $c_{i-1}$ and ends at $c_i$. 
Created arcs are marked by broken lines in Figure.
Clearly, $G_1 \subset \Gamma_1$, $G_i  \subset \Gamma_{i-1} \cup \Gamma_i$ for $i=2, \dots ,n-1$, and $G_n \subset \Gamma_{n-1}$.
We have constructed the pair of curves $F^* = F_1 \cup \dots \cup F_{n-1}$ and $G^* = G_1 \cup \dots \cup G_n$.

Denote by $U_\varepsilon$ the closure of the union of all arms of the butterflies of $R$ and of all pairs of curves of the form $F^*$ and $G^*$.
We see that $U_\varepsilon$ is obtained from ${\rm bd} (R)$ by exchanging all pairs of opposite curves $F$ and $G$ by the constructed pairs of curves $F^*$ and $G^*$.

We define $R_\varepsilon$ as the set bounded by $U_\varepsilon$.

\vskip0.2cm
Part B. We define a sequence of sets $R\sp j \subset S^2$ and show that they are convex bodies.

\vskip0.1cm
We define $R_0, R_1, \dots$ by induction. 

Put $R\sp 0 = R$.
Clearly it is a convex body.
 
Assume that $R\sp {j-1}$ is a convex body, where $j \geq 1$.
We get the boundary of $R\sp j$ by exchanging a pair of opposite curves $F, G$, if any remains, from $\bd (R\sp {j-1})$ into a pair $F\sp *, G\sp *$ as in Part A.

At every boundary point $p$, the set $R\sp j$ is supported by a hemisphere.  
If $p \in F\sp * \cup G\sp *$ this follows from the 
construction of $F\sp *$ and  $G\sp *$.
If $p \in {\rm bd} (R\sp j)$ does not belong to any curves $F^*$ or $G^*$, then in the part of the supporting hemisphere take this supporting $R\sp {j-1}$.
From Lemma 1 we conclude that $R\sp j$ is convex.
Clearly, $R\sp j$ is a convex body.

\vskip0.15cm
Part C. We show that $R_\varepsilon$ is a convex body and $\Delta (R_\varepsilon) = \Delta (R)$.

\vskip0.1cm
If after a finite number of steps in Part B we get $R_\varepsilon$, we see that it is a convex body. 

In the opposite case, $R_\varepsilon = {\rm lim}_{j\to \infty}R\sp j$. 
By Claim 1 we see that $R_\varepsilon$ is also a convex body.

From the construction of $R^j$ we see that for every supporting hemisphere $K$ of $R^j$ at every point of $F^* \cup G^*$ we have $\width_K = \Delta (R)$. 
The remaining part of $\bd(R^j)$ does not differ from $\bd(R^{j-1})$. 
So by induction we get $\Delta (R^j) = \Delta (R)$.
Thus if after a finite number of steps we obtain $R_\varepsilon$, it has thickness $\Delta(R)$.
By Lemma 4 of \cite{L2},
if $R_\varepsilon = {\rm lim}_{j\to \infty}R\sp j$, then the same is true.

\vskip0.15cm
Part D. Let us prove that $R_\varepsilon$ is a reduced body.

\vskip0.1cm
We should show that for any convex body $Z \subset R_\varepsilon$ different from $R_\varepsilon$ the inequality $\Delta (Z)  < \Delta (R_\varepsilon)$ holds true.
The body $Z$ does not contain an extreme point $e$ of $R_\varepsilon$ as it follows from the spherical analog of the Krein-Millman theorem formulated on p. 565 of \cite{L2}. 
Consequently, $Z$ is disjoint with an open disk $D$ centered at $e$. 
Now exactly as in Part 4 of \cite{L1} we show that $\Delta (Z)  < \Delta (R_\varepsilon)$, so that $R_\varepsilon$ is a reduced body.

\vskip0.15cm
Part E. We show that every $F^*$ and every $G^*$ are in the union of some triangles. 

\vskip0.1cm
Take a pair of curves $F^*$ and $G^*$ constructed in Part A.  
The bounding semicircle of the first lune supporting $R_\varepsilon$ at $f_i$ is denoted by $K_i$, and that at $g_i$ by $L_i$ (again see Figure).

For $i \in \{1, \dots , n-1 \}$, by $k_i$ denote the point of intersection of $K_i$ with $K_{i+1}$ (if they are subsets of a great circle, take $k_i$ as the center of $f_if_{i+1}$).
By $l_i$ denote the point of intersection of $L_i$ with $L_{i+1}$  
(if they are subsets of a great circle, take $l_i$ as the center of $g_ig_{i+1}$).

For every $c \in G_i$ take the lune supporting $R_\varepsilon$ such that $cf_i$ is the thickness chord and denote by $T(c)$ the bounding semicircle of this lune through $f_i$. 
In particular, $T(g_i) = K_i$.
When we move $c \in G_i$ counterclockwise from $g_i$ to $c_i$, the lune and thus also $T(c)$ ``rotate" counterclockwise.  
So since the distance between $c_i$ and any point of $T(c_i)$ is at least $\Delta(R)$ we see that the distance from $c_i$ to every point of the arc $f_ik_i$ is at least $\Delta(R)$. 
Analogously, the distance from $c_i$ to any point of $f_{i+1}k_i$ is at least $\Delta(R)$.
So since every point of $F_i$ is at the distance $\Delta(R)$ from $c_i$, we get $F_i \subset f_ik_if_{i+1}$.
Thus $F^*$ is contained in the union of triangles $f_ik_if_{i+1}$, where $i=1, \dots , n$.
Similarly, $G^*$ is in the union of triangles  $g_il_ig_{i+1}$, where $i=1, \dots , n$.

\vskip0.2cm
Part F. We show that the Hausdorff distance between $R$ and $R_\varepsilon$ is at most $\varepsilon$. 

\vskip0.1cm
Denote by $P$ the closure of the convex hull of all points $f_i$ and $g_i$ and of endpoints of all arms of butterflies of $R$. 
Denote by $Q$ the union of $P$ and all triangles $f_ik_if_{i+1}$ and $g_il_ig_{i+1}$.

Part E implies inclusions $P \subset R \subset Q \hskip0.3cm  {\rm and} \hskip0.3cm P \subset R_\varepsilon \subset Q$. 
So in order to estimate the Hausdorff distance between $R$ and $R_\varepsilon$ by $\varepsilon$ it is sufficient to estimate the Hausdorff distance between $P$ and $Q$ by $\varepsilon$.  
Since $P \subset Q$, the Hausdorff distance between them is $\sup_{q \in Q} \inf_{p \in P} |pq|$.
Hence it is  sufficient to show that all (i.e., for all pairs $F$, $G$ and all $i$) distances between $k_i$ and $f_if_{i+1}$, and also between $l_i$ and $g_ig_{i+1}$, are at most $\varepsilon$. 

Since the sum of angles of any quadrangle $o_if_ik_if_{i+1}$ is over $2\pi$, from the assumption in Part A on the angle between two successive chords below $\frac{\pi}{2}$, we get $\angle f_io_if_{i+1} < \frac{\pi}{2}$. 
So $\angle f_ik_if_{i+1} > 2\pi - 3 \cdot \frac{\pi}{2} = \frac{\pi}{2}$.
Hence $\tan \frac{1}{2}\angle f_ik_if_{i+1} > 1$. 
By Lemma 2 the distance from $k_i$ to $f_if_{i+1}$ is at most 

\vskip-0.1cm

$${\rm arc sin}\frac{\tan \frac{1}{2} |f_if_{i+1}|}{{\tan \frac{1}{2}  \angle f_ik_if_{i+1}}} <
{\rm arc sin} \;{\tan \frac{1}{2} |f_if_{i+1}|}.$$

By $|f_if_{i+1}| \leq \rho_\varepsilon$ (see Part A) and $\rho_\varepsilon = 2 \cdot {\rm arc tan} (\sin \varepsilon)$ we get that it is at most $\rm{arc sin} \tan \frac{\rho_\varepsilon}{2} = \rm{arc sin} \tan \frac{2 \cdot {\rm arc tan (\sin \varepsilon)}}{2} = \varepsilon$. 
Analogously, the distance between $l_i$ and $g_ig_{i+1}$ is at most $\varepsilon$.
We see that the Hausdorff distance between $P$ and $Q$, and thus between $R$ and $R_\varepsilon$, is at most $\varepsilon$.
\end{proof}

\vskip-0.1cm
If $R$ is a body of constant width, from Parts A and B of this proof we see that $R_\varepsilon$ (so $R^1$ in this case) is also a body of constant width. 
So we get the following corollary.

\vskip0.2cm
\noindent
{\bf Corollary.}
For every body $W \subset S^2$ of constant width and for arbitrary $\varepsilon >0$ there exists a body $W_\varepsilon \subset S^2$ of constant width $\Delta(W_\varepsilon) = \Delta (W)$ whose boundary consists only of arcs of circles of radius $\Delta(W)$, such that the Hausdorff distance between $W$ and $W_\varepsilon$ is at most $\varepsilon$.

\vskip0.2cm
Applications for Barbier's theorem in Part 4 of \cite {L1} do not hold true here for the sphere. 

Let us correct misprints in this Part 4 on p. 873: three times $2\pi$ should be changed into $\pi$.

\baselineskip 3pt

\vskip0.15cm
\noindent
Marek Lassak

\vskip0.05cm
\noindent
University of Science and Technology

\noindent
al. Kaliskiego 7, Bydgoszcz 85-796, Poland

\vskip0.1cm
\noindent
e-mail: lassak@utp.edu.pl


\baselineskip 6pt

\begin{thebibliography}{10}

\bibitem {Bl} W. Blaschke: Konvexe Bereiche gegebener konstanter Breite und kleinsten Inhalts, Math. Ann, 76 (1915) 504--513.

\bibitem {BF} T. Bonnesen, W. Fenchel: Theorie der konvexen K\"orper, Springer, Berlin et al., 1934, (Engl. transl. Theory of Convex Bodies, BCS Associated, Moscow, Idaho USA, 1987).

\bibitem {E} H. G. Eggleston: Convexity, Cambridge University Press, 1958.

\bibitem{FaLa} E. Fabi\'nska, M. Lassak: Reduced bodies in normed spaces, Isr. J. Math. 161 (2007) 75--88.

\bibitem{HA} N. N. Hai, P. T. Ann: A generalization of Blaschke's converegence theorem in metric space, J. Convex Anal. 20 1013--1024.  

\bibitem {L1} M. Lassak: Approximation of bodies of constant width and reduced bodies in a normed plane, J. Convex Anal. 19 (2012) 865--874.

\bibitem {L2} M. Lassak: Width of spherical convex bodies, Aeq. Math. 89 (2015), 555--567.

\bibitem{LasMart1} M. Lassak, H. Martini: Reduced convex bodies in Euclidean space - a survey, Expo. Math. 29 (2011) 204--21.

\bibitem{LasMart2} M. Lassak, H. Martini: Reduced convex bodies in finite-dimensional normed spaces -- a survey, Results Math. \textbf{66} (2014) 405--426. 

\bibitem {LaMu} M. Lassak, M. Musielak: Reduced spherical convex bodies, Bull. Pol. Acad. Sci., Math. 66 (2018) 87--97.

\end{thebibliography}
\end{document}